\documentclass[12pt]{amsart}

\usepackage{color}
\usepackage[centertags]{amsmath}
\usepackage[margin=1in]{geometry}
\usepackage{amsfonts}
\usepackage{amsthm}
\usepackage{newlfont}
\usepackage{verbatim}
\usepackage{amscd}
\usepackage{amsgen}
\usepackage[active]{srcltx}
\usepackage[all]{xy}
\usepackage{amssymb}
\usepackage{amssymb,amsmath}
\usepackage{amscd}
\usepackage{enumerate}
\usepackage{color}
\usepackage{xypic}
\usepackage{graphicx}
\usepackage[normalem]{ulem}

\newlength{\defbaselineskip} 
\theoremstyle{plain}
\newtheorem{thm}{Theorem}[section]
\newtheorem{cor}[thm]{Corollary}
\newtheorem{con}[thm]{Conjecture}
\newtheorem{lema}[thm]{Lemma}
\newtheorem{obs}[thm]{Proposition}

\theoremstyle{definition}
\newtheorem{df}[thm]{Definition}
\newtheorem{exm}[thm]{Example}
\newtheorem{rem}[thm]{Remark}

\newtheorem{fact}[thm]{Fact}

\newtheorem{pr}{Algorithm}

\theoremstyle{definition} 
\theoremstyle{definition} \newtheorem{ex}{Example}[section] %

 \numberwithin{equation}{section}
\def\p{\mathbb{P}}

\def\z{\mathbb{Z}}

\def\Z{\mathbb{Z}}

\def\c{\mathbb{C}}
\def\C{\mathbb{C}}

 \DeclareMathOperator{\Proj}{Proj}

\DeclareMathOperator{\codim}{codim} \DeclareMathOperator{\h}{ht}
\DeclareMathOperator{\Spec}{Spec}
\def\p{\mathbb{P}}
\def\a{\mathbb{A}}
\def\ob{\begin{obs}}
\def\kob{\end{obs}}
\def\dow{\begin{proof}}
\def\kdow{\end{proof}}

\def\tw{\begin{thm}}
\def\ktw{\end{thm}}
\def\hip{\begin{con}}
\def\khip{\end{con}}
\def\lem{\begin{lema}}
\def\klem{\end{lema}}
\def\ex{\begin{exm}}
\def\prog{\begin{pr}}
\def\kprog{\end{pr}}
\def\wn{\begin{cor}}
\def\kwn{\end{cor}}
\def\uwa{\begin{rem}}
\def\kuwa{\end{rem}}
\def\kex{\end{exm}}
\def\dfi{\begin{df}}
\def\kdfi{\end{df}}
 \definecolor{zielony}{rgb}{0.5, 0.9, 0.1}
 \definecolor{czerwony}{rgb}{0.9, 0.2, 0.1}
 \definecolor{niebieski}{rgb}{0.3, 0.1, 0.9}

\xyoption{line}
\def\fa{\begin{fact}}
\def\kfa{\end{fact}}

\definecolor{pink}{rgb}{1,0,1}

\newcommand{\A}[4]{\mathrm{A}^{#1,#2}_{#3,#4}}
\newcommand{\B}[4]{\mathrm{B}^{#1,#2}_{#3,#4}}
\newcommand{\M}{\mathrm{M}}
\newcommand{\g}{\mathfrak{g}}
\renewcommand{\h}{\mathfrak{h}}
\newcommand{\e}{\mathfrak{e}}
\renewcommand{\ll}{\mathfrak{l}}
\renewcommand{\1}{\mathbf{1}}
\newcommand{\m}{\mathbf{m}}
\newcommand{\aj}{\mathbf{a}}
\newcommand{\join}[2]{#1 \star #2}
\newcommand{\0}{\mathbf{0}}
\renewcommand{\AA}{\mathcal{A}}
\newcommand{\Ad}{\mathrm{adm}}
\newcommand{\X}{\mathrm{X}}

 \newif\ifprivate
\privatetrue

\def\???{\ifprivate {\bf {???}} \marginpar{{\Huge {\bf ?}}}
\else \fi}

\title{Local description of phylogenetic group-based models}

\author{Marta Casanellas}
\author{Jes\'us Fern\'andez-S\'anchez}
\author{Mateusz Micha{\l}ek}

\address{Departament de Matem\`atica Aplicada I \\ Universitat Polit\`ecnica de Catalunya \\
Av. Diagonal 647, 08028-Barcelona, Spain.}
 \email{marta.casanellas@upc.edu}
 \email{jesus.fernandez.sanchez@upc.edu}

 \address{Polish Academy of Sciences\\ ul.\'Sniadeckich 8, 00956 Warsaw, Poland} 
\email{mateusz.michalek@ujf-grenoble.fr} 

  \thanks{M. Casanellas and J. Fenr\'andez-S\'anchez are partially supported by
Spanish government MTM2012-38122-C03-01/FEDER and Generalitat de Catalunya 2009SGR1284. %
M. Micha\l ek was supported by Polish National Science Centre grant number DEC-2012/05/D/ST1/01063.
\newline{AMS 2000 subject classification 92D15;14H10;60J20}}

\address{}
 \email{}

\begin{document}
\begin{abstract}
Motivated by phylogenetics, our aim is to obtain a system of equations that define a phylogenetic variety on an open set containing the biologically meaningful points.
In this paper we consider phylogenetic varieties defined via group-based models. For any finite
abelian group $G$, we provide an explicit construction of $\codim X$ phylogenetic invariants (polynomial equations) of degree at most $|G|$  that define the variety $X$ on a Zariski open set $U$. The set $U$ contains all biologically meaningful points when $G$ is the group of the Kimura 3-parameter model. 
In particular, our main result confirms \cite[Conjecture 7.9]{jaPhD} and, on the set $U$, Conjectures 29 and 30 of \cite{SS}.
\end{abstract}
\maketitle
\section{Introduction}

As already devised  in the title of an essay  by J.E. Cohen (``Mathematics Is Biology's Next Microscope,
Only Better; Biology Is Mathematics' Next
Physics, Only Better" \cite{Cohen}), biology has lead to very interesting new problems in mathematics. In this paper we deal with algebraic varieties derived from phylogenetics, which were firstly introduced by Allman and Rhodes \cite{Allman2003}, \cite{Allman2004}. In phylogenetics, statistical models of evolution of nucleotides are proposed so that DNA sequences from currently living species are considered to have evolved from a common ancestor's sequence following a Markov process along a tree $T$. The living species are represented at the leaves of the tree, the interior nodes represent ancestral sequences, and the main goal in phylogenetics is to reconstruct the ancestral relationships among the current species. Roughly speaking, the phylogenetic variety $X$ associated to a Markov model and a tree $T$ is the smallest algebraic variety that contains the set of joint distributions of nucleotides at the leaves of the tree (see the introductory papers \cite{CasanellasNewsEMS} and \cite{AllmanRhodes_chapter4}). Its interest in biology lies in the fact that, no matter what the statistical parameters are, the (theoretical) joint distribution of nucleotides of the current species will be represented by a point in this phylogenetic variety. For this reason, the elements of the ideal of $X$ are known as \textit{phylogenetic invariants}. Knowing a system of generators of the ideal of $X$  would allow
doing phylogenetic inference without having to estimate the statistical parameters \cite{CFS_MBE}, which is always a tedious task.

Constructing a minimal system of generators of the ideal of phylogenetic invariants is hard and remains an open problem in most cases (for example, for the most general Markov model).  Apart from theoretical difficulties, its cardinal is huge (with respect to the number of leaves of the tree). On the other hand, a complete system of generators might have no biological interest, because the set of probability distributions forms
only a (real, semialgebraic) subset of the phylogenetic variety. Generalizing  some ideas of \cite{CFS} we propose a different approach: we construct a minimal system of $\codim X$ phylogenetic invariants that are sufficient to define $X$ on a Zariski open set containing the biological relevant points. We do this for certain phylogenetic varieties defined via the action of a finite abelian group $G$. These varieties turn out to be toric and comprise the phylogenetic varieties of two well known models in biology: Kimura 3-parameter model when $G=\z_2\times\z_2$, and the Felsenstein-Neyman model when $G=\z_2$. However, most of the well known models (such as Jukes-Cantor, Kimura 2-parameters, or the general Markov model) do not fit this description. A forthcoming paper involving different techniques will be devoted to these remaining models.

A toric variety $X$, equivariantly embedded in a projective space $\p$,
has a naturally distinguished open subset $U_X$: the orbit of the torus action.
When $X$ is the projective space $\p$, the corresponding open set $U:=U_{\p}$ is just the locus of points with all coordinates different from zero. In any case,
the variety $U_X$ is isomorphic to an algebraic torus (in particular it is smooth) and $U_X=X\cap U$. As it was observed in \cite{CFS} for the Kimura 3-parameter model, the biologically meaningful points of $X$ belong to $U_X$. It is thus well-justified to ask for the description of $U_X\subset U$.
The variety $U_X$ is in fact a complete intersection in $U$, hence it can be described by $\codim X$ phylogenetic invariants. We provide an explicit description of $U_X$ consistent with the following conjecture.

\hip\cite[Conjecture 29]{SS}\label{stop}
For any abelian group $G$ and any tree $T$, the ideal $I(X)$ of the associated phylogenetic variety is generated in degree at most $|G|$.
\khip
The conjecture is open, apart from the case $G=\z_2$ \cite{SS, Sonja}. The conjecture was stated separately for the Kimura 3-parameter model corresponding to the group $\z_2\times\z_2$ \cite[Conjecture 30]{SS}. In this case, it was proved first on the open subset $U$ \cite{JaAdvgeom} and later on the scheme-theoretic level \cite{JaJCTA}. In this paper, we give an explicit construction of phylogenetic invariants of degree at most $|G|$ that define the variety $X$ on $U$ for any finite abelian group $G$, providing insights to the above conjecture. Our proof has several steps. Starting from star trees and  cyclic groups, we inductively extend the construction to arbitrary abelian groups and trees. Moreover, we give a positive answer to a conjecture stated in the last author's PhD thesis \cite[Conjecture 7.9]{jaPhD}:

\hip\cite[Conjecture 7.9]{jaPhD}
On the orbit $U$ the variety associated to a claw tree is an intersection of varieties associated to trees with nodes of strictly smaller valency.
\khip


The paper is organized as follows. In section 2 we collect the preliminary results needed in the sequel and we give a local description of the phylogenetic varieties under consideration as a quotient of a group action. This result is supplementary to the rest of the paper but we include it because it sorts out an error in \cite{CFS}. In section 3 we provide the explicit generators of degree $\leq |G|$ for $U_X$ when $T$ is a tripod tree (first for the case of cyclic groups and then for arbitrary groups). In section 4, we give a construction for the desired generators of trees obtained by joining two smaller trees whose generators are already known. The results of these two sections provide the desired generators for the ideal of $U_X$ on any \textit{trivalent} tree $T$ (that is, a tree whose interior nodes have valency $\leq3$). In section 5 we stick to the case of claw trees of any valency, which allows us to provide the generators for arbitrary (not only trivalent) trees. Finally, in section 6 we describe the general procedure to obtain the desired generators for any tree and any abelian group, according to the results proved in the paper.


\subsection*{Acknowledgements}
The last author would like to thank Centre de Recerca Matemŕtica (CRM), Institut de Matemŕtiques de la Universitat de Barcelona (IMUB), Universitat Polit\`ecnica de Catalunya, and in particular Rosa-Maria Mir\'o-Roig, for invitation and great working atmosphere.

\section{Group-based models}

\subsection{Preliminaries}

An interesting introduction to tree models and its applications in phylogenetics can be found in \cite[Section 1.4.4]{PS} and \cite{AllmanRhodes_chapter4}. 
For the more specific case of group-based models we refer to \cite{SS, jaPhD} and we state the main facts here.

Let $G$ be a finite abelian group and $T$ a tree directed from a node $r$ that will be called the \textit{root}.
Let $E$, $L$ and $N$ be respectively the set of edges, leaves and interior nodes of the tree $T$.
Denote $\g:=|G|$, $\e:=|E|$ and $\ll:=|L|$, where $|\cdot|$ means cardinality. We assume the leaves of the tree are labelled so we have a bijection between $L$ and the set $\{1,2,\ldots,\ll\}$ and, in particular, an ordering on $L$. We will use additive notation for the operation in $G$.


\uwa
We assume that all the edges are directed from the root to simplify the notation. In fact, the orientation of edges of $T$ can be arbitrary -- all defined objects would be isomorphic \cite[Remark 2.4]{JaJalg}.
\kuwa

%

\dfi[group-based flow \cite{DBM,BW}]
A \textit{group-based flow} (or briefly a \textit{flow}) is a function $f:E\rightarrow G$ such that for each node $n$, we have $\sum_{e_i\in I_n}f(e_i)=\sum_{e_j\in O_n}f(e_i)$, where $I_n$ and $O_n$ are respectively the sets of edges incoming to $n$ and outgoing from $n$.
If the edges of the tree have been given an order $e_1, \dots,e_{\e}$, then a group-based flow $f$ will also be denoted by its values as $[f(e_1),\dots,f(e_{\e})].$
\kdfi

Notice that group-based flows form a group (isomorphic to $G^{\ll-1}$) with the natural addition operation -- cf. discussion on bijection between networks and sockets in \cite{JaJalg}.



\uwa\label{rem:leaves_flows}\rm
We want to point out that a group-based flow is completely determined by the values it associates to the pendant edges of $T$ (that is, if $f$ and $f'$ are two group-based flows such that $f(e_k)=f'(e_k)$ for all pendant edges $e_k$, then $f=f'$). Actually, given $g_1,\dots,g_l \in G$, there exists a flow $f$ that assigns $g_i$ at the pendant edge of leaf $i$ if and only if $\sum_i g_i=0$.
\kuwa

\ex\label{ex:flow}
Let us consider the group $G=\z_3$ and the following tree:%
\begin{figure}[h]
\begin{center}
 \includegraphics[scale=0.8]{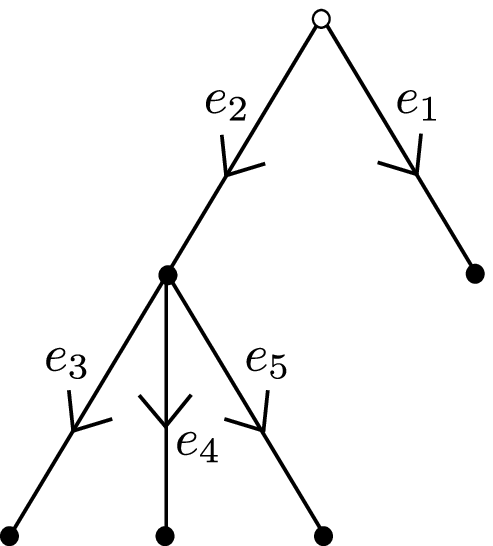}
\end{center}
\end{figure}
\newline
An example of a group-based flow on the tree above is given by the association $e_1\rightarrow 2, e_2\rightarrow 1$, $e_3\rightarrow 1$, $e_4\rightarrow 2$, $e_5\rightarrow 1$.
\kex

\dfi[Lattices $M$, $\tilde M$, Polytope $P_{T,G}$, Variety $X_{T,G}$, \cite{JaJalg}]\label{P} 
For a fixed edge $e_0$, we define $M_{e_0}$ as the lattice with basis elements indexed by all pairs $(e_0,g)$ for  $g\in G$. The basis element indexed by $(e_0,g)$ is denoted by $b_{(e_0,g)}\in M_{e_0}$.
%
We define $M=\prod_{j\in E} M_j$
and we have $M\simeq \z^{|E|\cdot|G|}$.
To a group-based flow $f$ one can naturally associate an element $Q_f:=\sum_{j\in E} b_{(j,f(j))}\in M$.
%
We define $P_{T,G}$ to be the polytope with vertices $Q_f$ over all group-based flows $f$ (for fixed $T$ and $G$ we will omit the subscript). The dimension of $P$ is taken as the topological dimension. 
 We also define $\tilde M$ to be the sublattice of $M$ spanned by $P_{T,G}$.


The variety $X_{T,G}:=\Proj\c[P]$ is called the \textit{phylogenetic variety associated to} $T$ \textit{and} $G$ (the $\Proj$ is over the semigroup algebra on the monoid generated by $P$ \cite{JaJalg}).
We will be omitting subscripts, if it does not lead to confusion.
\kdfi

\uwa
We consider the variety $X$ with its equivariant toric embedding, corresponding to the polytope $P$. This is not the same as the biologically meaningful embedding, but it is isomorphic (cf. \cite{SS}).
\kuwa



In terms of algebras, if we write $R=\c[x_{f_i}]$ with variables $x_{f_i}$ corresponding to group-based flows $f_i$ and $S=\c[\{y_{(e,g)}\}_{e\in E, g \in G}]$, then the ideal of $X$ is the kernel of the map:
$$\begin{array}{rcl}
R & \longrightarrow & S \\
x_{f} &\mapsto & \prod_{e \in E}y_{(e,f(e))}
\end{array}
$$
The associated map of affine spaces is
the \textit{parametrization of} $X$ \textit{according to the group-based model on} $G$.
It is worth noticing that $X_T$ is independent of the orientation of the tree (or placement of the root).
More precisely, if $T$ is an undirected tree and we root it at two different interior nodes $r, r'$ giving rise to
directed trees $T_r$ and $T_{r'}$, then $X_{T_r}$ is isomorphic to $X_{T_{r'}}$ \cite[Remark 2.4]{JaJalg}.


\ex
Some of the  varieties defined above come from biological evolutionary models: if we take  $G=\z_2$, we recover the \textit{Felsenstein-Neyman model} (or the binary Jukes-Cantor model), and  the \textit{Kimura 3-parameter model} corresponds to the group $G=\z_2\times\z_2$, see \cite{SS, HendyPenny}. For these groups, a discrete Fourier change of coordinates translates the parametrization map above into the original parametrization map used in biology, where the parameters stand for the transition probabilities between nucleotides.
\kex

By basic toric geometry, the equations of $X$ (thus elements of the defining ideal) correspond to relations among vertices of $P$ (see \cite[Chapter 13]{Stks}, \cite{Fult}, \cite{Cox}). These, by construction, correspond to group-based flows. 
The ideal of the corresponding phylogenetic variety $X$ is generated by those binomials $\prod_{i\in I} x_{f_i}=\prod_{j\in J} x_{f_j}$, such that $\sum_{i\in I} Q_{f_i}=\sum_{j\in J} Q_{f_j}$. A degree $d$ monomial $x_{f_{i_1}}x_{f_{i_2}}\dots x_{f_{i_d}}$ in the indeterminates of $R$ can be encoded as a multiset $m$ of $d$ group-based flows $m=\{f_{i_1},f_{i_2},\dots, f_{i_d}\}.$ Then each degree $d$  binomial in the ideal of $X$ is encoded as a relation between a pair of multisets $m=\{f_1, \dots, f_d\}$, $m'=\{f'_1, \dots, f'_d\}$  of $d$  flows each.
%
If $e_0$ is an edge of $T$, we denote by $\pi_{e_0}(m)$  the multiset $\{f_1(e_0), \dots, f_d(e_0)\}$.
 Then the multisets $m, m'$ correspond to a relation of degree $d$ among  flows (equivalently, to a phylogenetic invariant of degree $d$) if and only if the multisets $\pi_{e_0}(m)$ and $\pi_{e_0}(m')$ are equal for each edge $e_0\in E$. We denote this relation by $m \equiv_T m'$, $m\equiv_X m'$, or $m\equiv m'$ if the variety and the tree are understood from the context.
Then, we can write
\begin{eqnarray*}
 \{f_1,\ldots, f_d\}\equiv \{f'_1,\ldots, f'_d\} \quad \Leftrightarrow \quad  \sum_i Q_{f_i}=\sum_i Q_{f'_i}.
\end{eqnarray*}



\begin{figure}[h]
\begin{center}
 \includegraphics[scale=0.8]{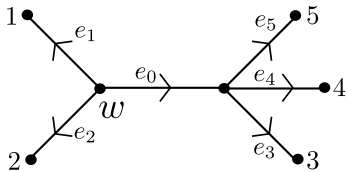}
\caption{\label{5-leaf tree} A directed 5-leaved tree.}
\end{center}
\end{figure}

\ex \label{example_5leaves}
Consider the binary Jukes-Cantor model (that is, $G=\z_2$) on the $T$ of Figure~\ref{5-leaf tree}.

In this case, group-based flows are represented by a sequence of $\e=6$ assignments of elements of $\z_2$ to the edges $e_0,\dots,e_5$ (in this order).
An example of a relation $m\equiv_T m'$ is given by the following pair of multisets, each containing two group-based flows:
\begin{enumerate}
\item $m=\{[1,1,0,1,0,0],[1,0,1,0,1,0]\}$,
\item $m'=\{[1,0,1,1,0,0],[1,1,0,0,1,0]\}$.
\end{enumerate}
The corresponding phylogenetic invariant is $x_{[1,1,0,1,0,0]}x_{[1,0,1,0,1,0]}=x_{[1,0,1,1,0,0]}x_{[1,1,0,0,1,0]}$.
\kex

\ex[Edge invariants]\label{edgeinvariants} The previous example is a special case of a general construction of edge invariants \cite{Pachter2004,AllRhMarkov}.
In our setting, edge invariants can be constructed as follows.
Fix an internal edge $e_0$ in a tree $T$, and decompose $T$ as a join of two trees $T_1$ and $T_2$
with only one common edge $e_0$ (this construction will be often used in Section \ref{gluetrees}).
Note that each group-based flow $f$ on $T$ decomposes exactly into two group-based flows $f_1,f_2$ respectively on $T_1$ and $T_2$, that assign the same element to $e_0$.
We will denote this by $f=f_1 \star f_2$.
Let us fix two group-based flows $f$ and $f'$ on $T$ that assign the same element to $e_0$ and decompose as $f=f_1 \star f_2$ and $f'=f_1'\star f_2'$.
Then, the following relation is called an \textit{edge invariant} associated to $e_0$ and is a quadratic phylogenetic invariant for the tree $T$:
\[\{f_1\star f_2, f_1' \star f_2'\} \equiv \{f_1 \star f_2', f_1'\star f_2\}.\]
Edge invariants can be defined on a more general setting (that is, for a broader class of evolutionary models) and correspond to  minors derived from rank conditions on  certain matrices associated to edges called \emph{flattenings} (see \cite{DK} and \cite{CFS3}).
\kex

\begin{figure}
\begin{center}
 \includegraphics[scale=0.8]{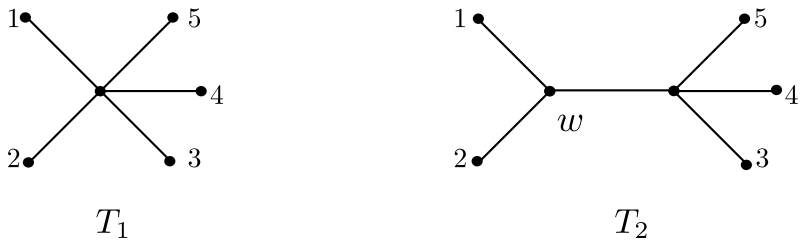}
\caption{\label{fig:2trees} $T_1\leq T_2$: $T_1$ is obtained by contraction of the interior edge of $T_2$. }
\end{center}
\end{figure}

\ex[Edge contraction]\label{subtree}
Given two undirected trees, we write  $T_1 \leq T_2$ if $T_1$ can be
obtained form $T_2$ by contraction of interior edges (that is, identifying both nodes of certain interior edges of $T_2$). For example, in figure \ref{fig:2trees} we have $T_1\leq T_2$.

If $T_1\leq T_2$, then the variety $X_{T_1,G}$ is contained in $X_{T_2,G}$ for any group $G$ (cf. for example \cite[Prop 3.9]{JaAdvgeom}) and the reverse inclusion holds for the corresponding ideals. That is, phylogenetic invariants of $T_2$ are also phylogenetic invariants of $T_1$. Note that phylogenetic invariants on $T_i$ have coordinates labelled by group-based flows on $T_i$. If we give to $T_1$ an orientation induced from an orientation on $T_2$, we can naturally associate to each flow on $T_2$ its restriction to $T_1$. For example, for the trees in figure \ref{fig:2trees}, we root $T_2$ at $w$ (obtaining the tree of Example \ref{example_5leaves}) and root $T_1$ at the interior node. Then a group-based flow $f$ on $T_2$ 
restricts to the flow on $T_1$ that assigns $f(e_i)$ to each pendant edge $e_i$ in $T_1$. For instance, if we consider  $G=\z_2$ and the invariant represented by $m \equiv_{T_2} m'$ as in Example \ref{edgeinvariants}, then the corresponding relation on $T_1$ is
$$\{[1,0,1,0,0],[0,1,0,1,0]\} \equiv_{T_1} \{[0,1,1,0,0],[1,0,0,1,0]\}$$
(that is, we have removed the assignment to $e_0$).

\kex

\subsection{Local description as a quotient}


The affine cone $\hat X$ over the variety $X$ equals $\Spec \c[P]$. The inclusion $P\subset \tilde M$ induces the inclusion of algebras $\c[P]\hookrightarrow\c[\tilde M]$. Geometrically, it gives the inclusion of the open set $\Spec \c[\tilde M]=U_{\hat X}\hookrightarrow \hat X$. Let $K_M\subset M$ be the positive quadrant in the lattice $M$.
The affine space of phylogenetic parameters equals $\mathbb{A}:=\Spec \c[K_M]$. The dominant map parameterizing $\hat X$ is induced by the inclusion $\c[P]\hookrightarrow\c[K_M]$. Again, one can restrict to dense torus orbits obtaining $U_{\mathbb{A}}:=\Spec \c[M]\rightarrow\Spec \c[\tilde M]$.

Our aim is to understand this map and its projectivization. As we will be considering projective varieties we introduce the following sublattices.
\dfi[Lattices $M_{0}$ and $\tilde M_{0}$]\label{MS0ME0}
We define $M_{0}$ as the sublattice of $M$ consisting of those points $p=\sum c_{(j,g)} b_{(j,g)}$ such that for each edge $j\in E$, the sum of coordinates $c_{(j,g)}$ over $g\in G$ equals $0$: $\sum_{g\in G} c_{(j,g)}=0$ for all $j\in E$. 

We define $\tilde M_{0}:=M_{0}\cap \tilde M$. This is the character lattice of the torus $U_X$ -- the locus of points of the projective toric variety $X$ with all coordinates different from zero \cite[Section 2.1]{Cox}. In other words, $U_X=X\cap U_\p\simeq \Spec\c[\tilde M_0].$
\kdfi

We have the following commutative diagrams:

\begin{eqnarray*}
\begin{array}{ccc}
\c[P] & \scalebox{1.5}{$\hookrightarrow$} & \c[K_M] \vspace*{-3mm}\\
 {\Large \rotatebox{270}{$\hookrightarrow$}} & & {\Large \rotatebox{270}{$\hookrightarrow$}} \\
 \c[\tilde M] & \scalebox{1.5}{$\hookrightarrow$} & \C[M] \\
  {\Large \rotatebox{90}{$\hookrightarrow$}} & & {\Large \rotatebox{90}{$\hookrightarrow$}} \\
 \c[\tilde{M}_0] & \scalebox{1.5}{$\hookrightarrow$} &  \C[M_0]
\end{array}
\qquad
\begin{array}{ccc}
\hat{X}=\mathrm{Spec}\, \c[P]& \longleftarrow & \mathbb{A}=\mathrm{Spec}\,\c[K_M]\\
 {\Large \rotatebox{90}{$\hookrightarrow$}} & & {\Large \rotatebox{90}{$\hookrightarrow$}} \\
 U_{\hat{X}}=\mathrm{Spec}\,\c[\tilde M] & \longleftarrow & U_{\mathbb{A}}=\mathrm{Spec}\,\C[M] \\
 \scalebox{1.5}{$\downarrow$} & & \scalebox{1.5}{$\downarrow$}\\
 U_X=\mathrm{Spec}\,\c[\tilde{M}_0]  & \longleftarrow & \mathrm{Spec}\,\C[M_0]
\end{array}
\end{eqnarray*}

\uwa\label{rem:localisation}\rm
If we consider the localization with respect to all the variables $x_{f_i}$ (that is, the Zariski open subset $U_X$ of $X$ isomorphic to a torus), we see that the variety is locally defined by Laurent binomials  $\left(\prod_{i\in I} x_{f_i}\right)\left(\prod_{j\in J} x_{f_j}\right)^{-1}-1$
such that $\sum_{i\in I} Q_i-\sum_{j\in J} Q_j=0$.
\kuwa


For a fixed node $n\in N$ we have an action of $G^* (\simeq G)$ on $\c[M]$ defined as follows: for $\chi\in G^*$ and given an indeterminate $y_{(j,g)}$ in $\c[M]$, $j\in E$, $g\in G$
we define :  
\begin{eqnarray*}
\chi \cdot_n y_{(j,g)}:=\left \{
\begin{array}{ll}
y_{(j,g)} & \mbox{if $j$ is not adjacent to $n$}, \\
\chi(g)\, y_{(j,g)} & \mbox{ if $j$ is incoming to $n$}, \\
\chi(g)^{-1}\, y_{(j,g)} & \mbox{ if $j$ is outgoing from $n$}.
\end{array}
\right .
\end{eqnarray*}
By the definition of the generators of $\tilde M$ we see that  $\c[\tilde M]$ is invariant by the action of $G^*$. 

Notice that for a fixed $n\in N$, the action $\cdot_n$ of $G^*$ restricts to $\c[M_0]$.

\ob\cite[Lemma 6.5 and Corollary 6.7]{jaPhD}\label{cor:codimension}
The following equality holds:
$$\c[\tilde M_0]=\c[M_0]^{(G^N)},$$
where $G^N\simeq (G^*)^N$ acts as above.
Moreover,
$$\dim P=\dim X=\dim U_X=\dim \tilde M_0=\dim M_0=(\g-1)\e,$$
hence
$$\codim X=\codim \hat X=\g^{\ll-1}-1-(\g-1)\e.$$
The projective parametrization map of the model, restricted to the dense torus orbits is a finite cover, given by a quotient of a finite group acting freely, where the cardinality of each fiber is
$$index (M_0:\tilde M_0)=|G^N|=\g^{|N|}.$$
\kob
The reader is referred to the Appendix for a proof of this result.
\uwa
The Proposition implies that each fiber of the projective parametrization map of the variety $X$ restricted to dense tori orbits has cardinality $\g^{|N|}$. Such fibers provide the possible parameters of the model.
\kuwa

Let us now prove that not only the cardinality of the fiber is constant. In fact, locally the parametrization map is described by the quotient of a free group action.
Indeed, consider the action of $(\c^*)^{\e}$ on $\c[M]$ given by
\begin{eqnarray*}
(\lambda_{e})_{e\in E} \cdot y_{(e_0,g)}=\lambda_{e_0} y_{(e_0,g)}.
\end{eqnarray*}
Notice that $M_0=M^{(\c^*)^{\e}}$. Consider a subtorus $(\c^*)^{\e-1}\subset(\c^*)^{\e}$ corresponding to points $(\lambda_{e})_{e\in E}$ such that $\prod \lambda_e=1$. The dense torus orbit in the affine space of parameters of the model is the spectrum of the algebra $\c[M]$, i.e. $U_{\mathbb{A}}=\mathrm{Spec}\, \c[M]$. By taking quotient in $U_{\mathbb{A}}$ additionally by $G^N\times(\c^*)^{\e-1}$ 
we obtain $U_{\hat{X}}$, or equivalently in the level of algebras:
\begin{eqnarray*}
\c[\tilde{M}]=\c[M]^{G^N\times(\c^*)^{\e-1}}.
\end{eqnarray*}
The group $G^N\times (\c^*)^{\e-1}$ acts also on $\c[K_M]$, the algebra of the whole parameter affine space $\c[y_{(e,g)}]_{e\in E,g\in G}$.
However, the quotient is not equal to the algebra of the affine cone over the variety $X$ representing the model (here the first two authors acknowledge an error in \cite[Theorem 3.6]{CFS} without further consequences in the quoted paper). 
Indeed, the algebra $\c[P]$ of the affine variety is invariant by the action of $G^N\times (\c^*)^{\e-1}$. 
However, the invariant monomials of $\c[K_M]$ correspond to all the monomials of $\tilde M$ that are in the positive quadrant of $M$.
Not all such monomials are generated by the polytope. For example, for the  Kimura $3$-parameter model (that is, $G=\z_2 \times \z_2$) the monomial $y_{(e_0,g)}^2\prod_{e_i\in E} y_{(e_i,0)}^2$, where $0\in G$ is the neutral element is invariant for any $g\in G$ and any distinguished edge $e_0$ (because $g+g=0$). This is not however the sum of any two vertices of the polytope associated to the variety. This is the reason why the quotient construction holds only locally in the Zariski topology.


\section{Complete intersection for the tripod}\label{sec:tripod}
By definition, the tripod is the tree $T$ with one interior node and three leaves. Assume the tripod is rooted at the interior node and leaves are labelled.
In this case, group-based flows can be identified with triples of group elements $[i,j,k]$ such that $i+j+k=0$. 
Consider a $\g\times \g$ matrix $\M=(m_{i,j})_{i,j}$ with columns and rows indexed by the group elements and only integral entries.
The entry at $(i,j)$ corresponds to the flow $[i,j,-i-j]$, and so, the matrix $\M$ can be identified with the Laurent monomial:
$$L(\M):=\prod_{i,j\in G} x^{m_{i,j}}_{[i,j,-i-j]}.$$

From now on, we will write $\mathcal{M}_{k}(\z)$ for the group of $k\times k$ matrices with integral entries under addition.
\lem\label{lem:admissible}
For a given matrix $\M\in \mathcal{M}_{\g}(\z)$, the Laurent binomial $L(\M)-1$ belongs to the localized ideal of the phylogenetic variety of $X_T$ if and only if:
\begin{enumerate}
\item each row sum in $\M$ equals zero,
\item each column sum in $\M$ equals zero,
\item for each $k\in G$, the sum of all entries $m_{i,j}$ with $i+j=k$ equals zero.
\end{enumerate}

\klem
\dow
Follows from Remark \ref{rem:localisation} -- the three conditions correspond to three different edges of the tripod.
\kdow

\dfi\label{def:admissible} \rm
Integral matrices satisfying the three conditions of Lemma  \ref{lem:admissible} will be called \emph{admissible for $G$}. That is, if for any $k\in G$ we define the set $S_k=\{(i,j)\in G\times G \mid i+j=k\}$, a matrix $\M$ is admissible (for $G$) if 
\begin{enumerate}
\item each row sum in $\M$ equals zero,
\item each column sum in $\M$ equals zero,
\item for each $k\in G$, $\sum_{(i,j)\in S_k}m_{i,j}=0$.
\end{enumerate}
Given any matrix (admissible or not), we call the sum of positive entries its \emph{degree}. It is equal to the degree of the associated binomial, obtained from $L(\M)-1$ after clearing denominators:
\[\prod_{m_{i,j}>0} x_{[i,j,-i-j]}^{m_{i,j}}-
\prod_{m_{i,j}<0} x_{[i,j,-i-j]}^{-m_{i,j}}.\]
\kdfi

\ex
Consider the group $G=\z_3$. The order of elements labelling the rows and columns is as follows: $0,1,2$. Consider the following, matrix with degree three:
\begin{eqnarray*}
\M=\left( \begin{array}{ccc}
0 & -1 & 1 \\
1 & 0 & -1 \\
-1 & 1 & 0 \\
\end{array} \right).
\end{eqnarray*}
It corresponds to the monomial:
\begin{eqnarray*}
 L(\M)=
 x_{[0,1,2]}x_{[1,2,0]}x_{[2,0,1]}x_{[0,2,1]}^{-1}x_{[1,0,2]}^{-1}x_{[2,1,0]}^{-1},
\end{eqnarray*}
and to the binomial:
$$x_{[0,1,2]}x_{[1,2,0]}x_{[2,0,1]}-x_{[0,2,1]}x_{[1,0,2]}x_{[2,1,0]}.$$
\kex

\dfi[$\Ad(G)$]
Any integral combination of admissible matrices gives an admissible matrix, so that admissible $\g\times \g$ matrices form a $\z$-submodule of $\mathcal{M}_{\g}(\z)$. We will denote it by $\Ad(G)$.
\kdfi

\uwa \rm
Admissible matrices of a cyclic group have an interpretation in terms of ``magic squares". They are differences of two magic squares with:
\begin{itemize}
\item[(i)] each row summing up to a fixed number,
\item[(ii)] each column summing up to a fixed number,
\item[(iii)] each generalized diagonal (that is $\g$ entries parallel to the diagonal) summing up to a fixed number.
\end{itemize}
\kuwa

\wn\label{cor:admis}
A set of Laurent binomials $L(\M_i)-1$, $i=1,\dots, m$, defines the variety $X$ in the Zariski open set $U$ if and only if $M_i$ generate $\Ad(G)$.
\kwn

\subsection{Cyclic group}

Consider $G=\z_\g$. In this case, the codimension of the variety $X$ associated to the tripod is $(\g-1)(\g-2)$. We shall explicitly construct $(\g-1)(\g-2)$ binomials of degree at most $\g$ that define $X$ in the Zariski open set $U$.
Each binomial will have the form $L(\M)-1$, where $L(\M)$ is the Laurent monomial corresponding to an admissible matrix $\M$. By virtue of Corollary~\ref{cor:admis} we need to construct $(\g-1)(\g-2)$ admissible matrices for $G$ that form a $\z$-basis of $\Ad(G)$.

\uwa
The construction we give below works for any value of $\g$. However, for values of $\g$ that are not powers of a prime number it may be better to use the method of Section \ref{arbitgr} in order to obtain lower degree phylogenetic invariants.
\kuwa

We begin by  introducing elementary matrices in $\mathcal{M}_{\g}(\z)$ as follows: given $i,j \in G$, we write $E^{i}_{j}\in \mathcal{M}_{\g}(\z)$ for the matrix whose entries are all equal to zero, except for the entry in the $i$ row and $j$ column, which is equal to one.


\dfi[Matrix {$\A{i}{a}{j}{b}$} ] \rm
Given group elements $i,j,a,b \in G$ with  $i \neq a$, $j\neq b$, we define the matrix
\begin{eqnarray}\label{def:A}
\A{i}{a}{j}{b}=
 (E^{i}_{j} + E^{a}_{b})-(E^{i}_{b}+E^{a}_{j}),
 \end{eqnarray}
 that is, 
\begin{eqnarray*}
\qquad  j  \qquad \quad  b \qquad \qquad \\
 \A{i}{a}{j}{b}=\begin{array}{c} \\ i\\ \\a \\ \\ \end{array} \, \left (
 \begin{array}{ccccc}
  & \vdots & & \vdots & \\
  \cdots & 1 & \cdots & -1 & \dots \\
  & \vdots & & \vdots & \\
  \cdots & -1 &\cdots  & 1 & \cdots \\
  & \vdots & & \vdots &
   \end{array}
 \right )
\end{eqnarray*}
where all the non-explicit entries are zero.
\kdfi
This matrix has degree 2  and, although it is not admissible, it satisfies the first two conditions of admissibility. Our idea is to produce phylogenetic invariants, or equivalently admissible matrices, as sums of these matrices.

\ob\label{prop:cyclic}
For $G=\Z_\g$ there exists an integral basis of cardinality $(\g-1)(\g-2)$ of the $\Ad(G)$, where each matrix has degree at most $\g$.
\kob
\dow

First, we construct the candidate to basis.
We define the following set
\begin{eqnarray*}
K=\{(i,j)\in G\times G \mid i\neq 0 \textrm{ and } j\neq 0,1 \},
\end{eqnarray*}
which has cardinal $(\g-1)(\g-2).$
For each index $(i,j)\in K$ we shall define an admissible matrix $\X(i,j)$ that has $1$ on the entry corresponding to $(i,j)$ and $0$ on all the other entries indexed by  $K$:
\begin{eqnarray*}
&  \begin{array}{ccccc}
\; \; 0 &  1 & \quad \quad & j & \qquad
\end{array} & \\
 \X(i,j)= \begin{array}{c}
 0  \\
\; \\
i \\
\;  \\ \\
\end{array}
&
 \left [
 \begin{array}{cc|ccc}
 * & * &  * &  \ldots & * \\ \hline
 * & * &    & \vdots  &   \\
 \vdots & \vdots &   \ldots & 1  & \ldots  \\
 * & * &    & \vdots  &   \\
      \end{array}
 \right ]  &
\end{eqnarray*}

\begin{figure}
\begin{eqnarray*}
\X(2,4)=
 \left [
 \begin{array}{cc|ccccc}
0 & 1 &  0 &  0 & -1 & 0 & 0\\ \hline
 -1 & 1 &  & & \vdots & &  \\
 -1 & 0 &  &\cdots  & 1 & &  \\
 0 & 0 &  & &  & &  \\
 1 & -1 &  & &  & &  \\
 1 & -1 &  & &\ & &  \\
 0 & 0 &  & & & &  \\
      \end{array}
 \right ]
 \end{eqnarray*}
 \caption{Example of the matrix $\X(i,j)$ when $G=\z_7$,  $i=2$ and $j=4$. The entries that are not specified are zero.}
\end{figure}

We need to distinguish three cases:

\vspace*{2mm}

 \paragraph{\textbf{Case 1 ($i < \g/2$): }}
 We define $\X(i,j)$ as the following sum:
 \begin{eqnarray}\label{def_X}
  \A{i}{0}{j}{0}+\A{i-1}{j}{1}{0}+\A{i-2}{j+1}{1}{0}+\dots+\A{1}{i+j-2}{1}{0}+\A{0}{i+j-1}{1}{0}.
 \end{eqnarray}
This matrix has 1 at $(i,j)$ and zero at all other entries indexed by  $K.$ Moreover, it is admissible. Indeed, conditions (1) and (2) of Definition \ref{def:admissible} are satisfied by all the  summands above, so also for $\X(i,j)$.
%
On the other hand, using the definition of the $A$-matrices in \ref{def:A}, it is easy to see that $\X(i,j)$ can be decomposed as a sum of $i+1$ matrices as follows
\begin{eqnarray*}
 \X(i,j) = (E^{i}_{j}-E^{i+j-1}_{1}) + (E^{j}_{0}-E^{0}_{j})+ \sum_{a=0}^{i-1} (E^{a}_{1}-E^{a+1}_{0}) +\sum_{a=j}^{i+j-1}(E^{a}_{1}-E^{a+1}_{0}).
\end{eqnarray*}
Now, each difference between brackets satisfies the condition (3) and so, is an admissible matrix. It follows that $\X(i,j)$ is an admissible matrix.

As $\X(i,j)$ is a sum of $i+1$ matrices of degree $2$, its degree is less than or equal to $2(i+1)$. In fact, as the entry $(0,0)$ equals $1$ for the first matrix in the sum and $-1$ for the last matrix, the degree of $A$ is at most $2i+1.$ Our assumption in this case was $i <\g/2$, so that  $A$ is an admissible matrix of degree at most $\g$.

\vspace*{2mm}
\paragraph{\textbf{Case 2 ($i > \g/2$): }} We define $\X(i,j)$ as
\[\A{i}{0}{j}{0}+\A{j-1}{i}{1}{0}+\A{j-2}{i+1}{1}{0}+\dots+\A{j-(\g-i)}{i+(\g-i)-1}{1}{0}.\]
The same argument above proves that $\X(i,j)$ is an admissible matrix (taking into account that $j-\g+i=j+i$ and $i+\g-i=0$ in $\Z_\g$) of degree at most $2(\g-i+1)$. But as the entry $(i,0)$ appears once with $+1$ and once with $-1$, the degree of $\X(i,j)$ is actually less than or equal to $2(\g-i)+1$. Now we are assuming $i > \g/2$, so this degree is at most $\g$.

\vspace*{2mm}
\paragraph{\textbf{Case 3 ($i = \g/2$): }}
If $j \neq \g/2$ we proceed analogously to the cases above but exchanging the roles of $i$ and $j.$
Hence, the only case left is $\g=2k$ and $i=j=k$. In this case we define $\X(i,j)$ as in case 1. In $\X(i,j)$ we have the sum of $k+1$ degree $2$ matrices. However, each entry $(k,0)$ and $(0,0)$  appears twice with different signs, thus the matrix is of degree at most $\g$.

Now we consider the set $\mathcal{B}=\left\{\X(i,j)\mid (i,j)\in K \right\}$ and prove that it is indeed a $\Z$-basis for $\Ad(G)$.
Matrices in $\mathcal{B}$ are clearly linearly independent because they have one entry equal to 1 and all other entries labelled by $K$ equal to zero. Moreover, matrices in $\mathcal{B}$ generate any admissible matrix $\M$. Indeed, by subtracting an integral combination of the matrices in $\mathcal{B}$ we obtain an admissible matrix with all entries labelled by $K$ equal to zero. It is an easy observation that such a matrix is necessarily zero. Hence, $\M$ is an integral combination of the  matrices in $\mathcal{B}$.
\kdow
\uwa
Most of the phylogenetic invariants constructed above are of degree smaller than $\g$, however some of them may be exactly of degree $\g$. This complies with the conjecture of Sturmfels and Sullivant \cite{SS}.
\kuwa
\ex
For the group $\z_4$ the above construction gives the following $6$ matrices:
\footnotesize
\begin{eqnarray*}
 \X(1,2)=\left( \begin{array}{cc|cc}
0 & 1 & -1 & 0 \\ \hline
-1 & 0 & 1 & 0 \\
1 & -1 & 0 & 0 \\
0 & 0 & 0 & 0 \\
\end{array} \right), \qquad
\X(1,3)= \left( \begin{array}{cc|cc}
0 & 1 & 0 & -1 \\ \hline
-1 & 0 & 0 & 1 \\
0 & 0 & 0 & 0 \\
1 & -1 & 0 & 0 \\
\end{array} \right),\\
\X(2,2)= \left( \begin{array}{cc|cc}
1 & 0 & -1 & 0 \\ \hline
-1 & 1 & 0 & 0 \\
0 & 0 & 0 & 0 \\
0 & -1 & 1 & 0 \\
\end{array} \right),\qquad
\X(2,3)= \left( \begin{array}{cc|cc}
1 & 0 & 0 & -1 \\ \hline
0 & 0 & 0 & 0 \\
-1 & 1 & 0 & 0 \\
0 & -1 & 0 & 1 \\
\end{array} \right),\\
 \X(3,2)=\left( \begin{array}{cc|cc}
1 & 0 & 0 & -1 \\ \hline
-1 & 1 & 0 & 0 \\
-1 & 0 & 0 & 1 \\
1 & -1 & 0 & 0 \\
\end{array} \right),\qquad
 \X(3,3)=\left( \begin{array}{cc|cc}
0 & 1 & -1 & 0 \\ \hline
-1 & 1 & 0 & 0 \\
0 & -1 & 1 & 0 \\
1 & -1 & 0 & 0 \\
\end{array} \right).
\end{eqnarray*}
\normalsize
The matrices $\X(1,2)$ and $\X(1,3)$ have degree $3$ and correspond to Case 1;  the matrices $\X(2,2)$ and $\X(2,3)$ correspond to Case 2, and the last two $\X(3,2)$ and $\X(3,3)$ to Case 3. The rows and columns are labeled consecutively with group elements $0,1,2,3$.
\kex

\subsection{Arbitrary group}\label{arbitgr}

Consider the group $G\times H$, where $G$ and $H$ are \emph{arbitrary} abelian groups, and denote $\g=|G|$ and $\h=|H|$. Suppose that we already know generating sets for $\Ad(G)$ and $\Ad(H)$, which have cardinalities $(\g-1)(\g-2)$ and $(\h-1)(\h-2)$, respectively. We shall define a collection of $(\g\h-1)(\g\h-2)$ matrices, which will be a generating set for $\Ad(G\times H)$. Moreover, these admissible matrices either will have degree three or will come from admissible matrices of $G$ or $H$. In particular, the degree of the equations does not increase, apart from the fact that we are adding new cubics.

In this subsection we will use the following notation: elements of $G$ will be denoted by $i,j$ and always precede the comma when a couple $(\cdot, \cdot)\in G\times H$ occurs; elements of $H$ will be called $k,l$ and always follow the comma. We point out that, although we are using the same notation for the group elements as we did in the cyclic case, the groups $G$ and $H$ are not necessarily cyclic.

\dfi[matrix {$\B{i}{j}{k}{l}$}]
Given elements $i,j\in G$ and $k,l\in H$, with $j\neq 0$ and $k\neq 0$, we define the admissible matrix $\B{i}{j}{k}{l}$ for $G\times H$ by
\begin{eqnarray*}
 \B{i}{j}{k}{l}& = & \left(E^{(i,k)}_{(j,l)}+E^{(i+j,0)}_{(0,l)}+E^{(i,0)}_{(0,k+l)} \right)- \left (E^{(i+j,0)}_{(0,k+l)}+E^{(i,k)}_{(0,l)}+E^{(i,0)}_{(j,l)}\right),
\end{eqnarray*}
that is,
\begin{eqnarray*}
\quad  (0,l)  \qquad (0,k+l) \qquad (j,l)  \qquad \quad \\
 \B{i}{j}{k}{l}=
 \begin{array}{c} \\ (i,0) \\ \\(i,k) \\ \\ (i+j,0) \\  \\
 \end{array} \, \left (
 \begin{array}{ccccccc}
  & \vdots & & \vdots & & \vdots & \\
  \cdots & 0 & \cdots & 1 & \cdots & -1 & \cdots\\
  & \vdots & & \vdots & & \vdots & \\
  \cdots & -1 & \cdots & 0 & \cdots & 1 & \cdots\\
  & \vdots & & \vdots & & \vdots & \\
\cdots & 1 & \cdots & -1 & \cdots & 0 & \cdots\\
  & \vdots & & \vdots & & \vdots &
   \end{array}
 \right )
\end{eqnarray*}
where rows and columns have been ordered lexicografically and all non-explicit entries are zero.
%
It is a degree $3$ matrix, representing a cubic polynomial.
\kdfi

Let us consider the following sets of admissible matrices:
\begin{enumerate}
\item $(\g-1)^2(\h-1)^2$ matrices $\B{i}{j}{k}{l}$ for $i,j\neq 0$ and $k,l\neq 0$;
\item $(\h-1)^2(\g-1)$ matrices $\B{0}{j}{k}{l}$ for $j\neq 0$ and $k,l\neq 0$;
\item $(\g-1)^2(\h-1)$ matrices $\B{i}{j}{k}{0}$ for $i,j\neq 0$ and $k\neq 0$;
\item $(\h-1)^2(\g-1)$ transpositions of matrices of type (2);
\item $(\g-1)^2(\h-1)$ transpositions of matrices of type (3);
\item $(\g-1)(\h-1)$ matrices $\B{0}{j}{k}{0}$ for $j\neq 0$ and $k\neq 0$;
\item a $\z$-basis for $\Ad(G)$ embedded in $\mathcal{M}_{\g\h}(\z)$ by putting the elements of $H$ equal to zero (this produces $(\g-1)(\g-2)$ matrices);
\item a $\z$-basis for $\Ad(H)$ embedded in $\mathcal{M}_{\g\h}(\z)$ by putting the elements of $G$ equal to zero (this produces $(\h-1)(\h-2)$ matrices).
\end{enumerate}
Let us explain the construction of matrices of type $(7)$ and $(8)$. There is a canonical embedding $s_G:G\rightarrow G\times H$ defined by $s_G(i)=(i,0)$. By this embedding we can regard elements of $G$ as elements of $G\times H$. Hence, we can regard an admissible matrix for $G$ as a submatrix of an admissible matrix for $G\times H$, by putting all entries not indexed by group elements in the image of $s_G$ equal to zero. This gives matrices of type $(7)$ and, analogously, of type $(8)$.

All together we have defined a set $\mathcal{B}$ of $(\g\h-1)(\g\h-2)$ matrices. 
We prove below that any admissible matrix of $G\times H$ is an integral combination of matrices in $\mathcal{B}$.

\uwa
Note that the group of group-based flows acts on the coordinates of the ambient space $\mathbb{A}$. Consider a matrix $\B{0}{j}{k}{0}$ representing a simple cubic relation among flows (type (6) above). Any other relation $\B{i}{j}{k}{l}$ is obtained from it by the action of the group-based flow $[(i,0),(0,l),(-i,-l)]$. %
The action of the  group $\mathfrak{S}_3$ on the leaves of $T$ induces an action on the coordinates of the ambient space. Namely, if $\sigma\in \mathfrak{S}_3$ and $[f_1,f_2,f_3]$ is a group-based flow, define
\begin{eqnarray*}
 \sigma \cdot x_{[f_1,f_2,f_3]} =  x_{[f_{\sigma(1)},f_{\sigma(2)},f_{\sigma(3)}]}.
\end{eqnarray*}
Similarly, the transposition of the matrix $\B{i}{j}{k}{l}$ corresponds to the action of the transposition $\tau_{12}=(12) \in \mathfrak{S}_3$. Thus, although it may seem that the types (1)-(6) of matrices introduced above are complicated, they all come from the most simple type (6) matrices $\B{0}{j}{k}{0}$ under the action of a group equal to the semidirect product of $\mathfrak{S}_3$ and the group of group-based flows.
\kuwa
\ob
Each admissible matrix for the group $G\times H$ is an integral combination of admissible matrices of types $(1)-(8)$.
\kob
\dow
Consider an admissible matrix $\M$ of $G\times H$. We will reduce it to zero modulo the matrices presented above. First note that the matrices of type (1) have  a unique nonzero entry (equal to 1) indexed by $(i,k)(j,l)$ for $i,j\neq 0$ and $k,l\neq 0$. 
Thus, by subtracting an integral combination of matrices of type (1)  we can reduce all such entries to zero. We proceed analogously for entries indexed by $(0,k)(j,l)$, $(i,k)(j,0)$,$(i,k)(0,l)$, $(i,0)(j,l)$ for $i,j\neq 0$ and $k,l\neq 0$, using respectively matrices of type (2), (3), (4) and (5). Hence, we can assume that the only nonzero entry of the matrix $\M$ are indexed  by $(i,k)(j,l)$ where at least two of $i,j,k,l$ are neutral elements in the groups to which they belong. Entries indexed by $(0,k)(j,0)$ for $k\neq 0$, $j\neq 0$ can be reduced using matrices of type (6).

Notice that we did not reduce entries indexed by $(i,0)(0,l)$ or $(i,k)(0,0)$ nor $(0,0)(j,l)$. 
We claim that if $i,l\neq 0$, these entries are in fact 0. Indeed, fix a column indexed by $(i,l)$ with $i\neq 0$, $l\neq 0$. After the reduction process described above, we know that the only possible nonzero entry in this column is $(0,0)(j,l)$. By admissibility this entry must also be zero. The same holds for $(i,k)(0,0)$ by considering a row. Consider now an entry indexed by $(i,0)(0,l)$ for $i\neq 0$ and $l\neq 0$. The sum of these two indices equals $(i,l)$ but no sum of indices of any other nonzero entry in the matrix is equal to $(i,l)$ (all remaining entries have indices of type $(0,\cdot)(0,\cdot)$ or $(\cdot,0)(\cdot,0)$, which do not sum up to $(i,l)$ if $i,l\neq0$).
Thus, by admissibility, the entry indexed by $(i,0)(0,l)$ must be equal to zero.

Hence, we have reduced the matrix to a matrix $\M$ that has nonzero entries indexed only either by $(0,k)(0,l)$ or $(i,0)(j,0)$ for some (possibly equal to $0$) elements $i,j,k,l$.
It remains to show that such a matrix is a sum of the admissible matrices induced from $G$ or $H$. This will finish the proof, as such matrices, by assumption, are integral combinations of matrices of type $(7)$ and $(8)$.

Let $S_1$ be the subset of entries indexed by $(i,0)(j,0)$ for $i,j\in G$ and let $S_2$ be the subset of entries indexed by $(0,k)(0,l)$ for $k,l\in H$. The intersection $S_1\cap S_2$ contains precisely one entry indexed by $(0,0),(0,0)$, which we call $e$.
Let us define a matrix $\M_1$, which will be an admissible matrix induced from $G$, as follows. Each entry in $S_1$ different from $e$ is defined to be the same in $\M_1$ and $\M$. Moreover, all entries of $\M_1$ not in $S_1$ are set to zero. It remains to define the entry $e$. We define it so that the sum of the row indexed by $(0,0)$ in $\M_1$ is equal to zero. Let us notice that all other rows of $\M_1$ either coincide with $\M$ or have all entries equal to zero. Hence, the sum of all entries of $\M_1$ is equal to zero. For the same reason, all columns of $\M_1$ not indexed by $(0,0)$, have entries summing up to zero. Hence, so must the column indexed by $(0,0)$ and $\M_1$ satisfies the first two conditions of admissibility. We proceed to check the third condition. If $k\neq 0$ then all the entries with indices in $I_{(i,k)}$ for $\M_1$ are zero, and in particular sum up to zero. If $i\neq 0$ then all the entries indexed by elements of $I_{(i,0)}$ for $\M_1$ coincide with entries in $\M$, thus they sum up to zero. As the sum of all entries of $\M_1$ is zero, it follows that the sum of entries indexed by elements of $I_{(0,0)}$ equals zero. Thus $\M_1$ is admissible.
It immediately follows that $\M-\M_1$ is admissible and induced from $H$ (as all its nonzero entries are in $S_2$), which finishes the proof.
\kdow

Summing up, we have obtained the following result.

\begin{thm} For any finite abelian group $\z_{a_1}\times\dots\times\z_{a_k}$, the variety $X_T$ for the tripod $T$ is defined in $U$ by a complete intersection whose equations have degree at most $\max(3,a_i)$ and can be derived by successively applying the previous results.
\end{thm}

\begin{rem}
Note that for the 3-Kimura model we have constructed a set of phylogenetic invariants of degree $3$ that do not generate the whole ideal, but define the variety on an open set. This improves previous results from \cite{CFS} in the sense that invariants defining the variety on the open set given in that paper had degrees 3 and 4.
\end{rem}

\section{Complete intersection for joins of trees}\label{gluetrees}
\begin{figure}
\begin{center}
 \includegraphics[scale=0.3]{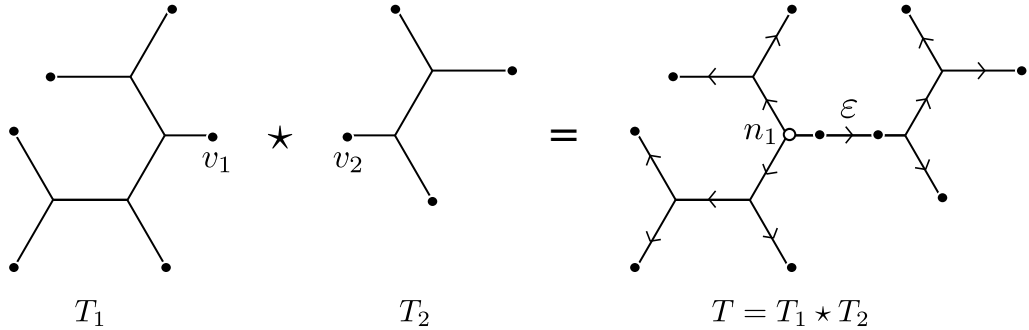}
 \label{gluing}
\caption{Gluing two trees $T_1$ and $T_2$ by the leaves $v_1$ and $v_2$. The resulting tree $T=T_1\star T_2$ is rooted at the node $n_1$, which is the closest node to $v_1$ in $T_1$, and oriented accordingly.}
\end{center}
\end{figure}

Let $G$ be an arbitrary finite abelian group of cardinality $\g$. Consider two (rooted or unrooted) trees $T_1$ and $T_2$, each with a distinguished leaf, say $v_1$ and $v_2$. We define the \emph{joined} tree $T=T_1\star T_2$ to be the tree obtained by identifying the leaves $v_1$, $v_2$ to an edge $\varepsilon$ (see Figure \ref{gluing}). It is well known how to find phylogenetic invariants for $T$ knowing them for $T_1$ and $T_2$ \cite{SS, SethTFP, JaJalg}. However, it is not clear how to describe the variety $X_T$ as a complete intersection in the Zariski open subset $U$, knowing such description for $X_{T_1}$ and $X_{T_2}$. This is the goal of this section.

The following notations and assumptions will be adopted throughout this section without further reference.
Write $\ll_i$ for the number of leaves of $T_i$ and $\e_i$ for the number of edges, so that $T=T_1\star T_2$ has $\ll:=\ll_1+\ll_2-2$ leaves and $\e:=\e_1+\e_2-1$ edges. We root $T$ at the node $n_1$ of $T_1$ closest to $v_1$ and orient $T$ from this root. The trees $T_1$ and $T_2$ will be given the orientation induced from this orientation on $T$ (see Figure \ref{gluing}). Moreover, for each tree $T_i$ choose a leaf $l_i$ different from $v_i$.


\dfi[$E_1(\cdot)$,$E_2(\cdot)$]\rm
Consider a group-based flow $f$ on $T_1$. There exists precisely one group-based flow $E_1(f)$ on $T$ that agrees with $f$ on $T_1$ and associates to all other leaves, apart from $l_2$, the neutral element of $G$. Indeed, take $E_1(f)$ equal to $f(n_1\rightarrow v_1)$ at all edges that appear in the shortest path from $n_1$ to $l_2$, and equal to the neutral element at the other edges of $T_2$. We call $E_1(f)$ \textit{the extension of $f$ to $T$ (relative to $l_2$)} . Analogously, for a group-based flow on $T_2$ we define the extension $E_2(f)$ (relative to $l_1$).
\kdfi



\begin{figure}
\begin{center}
 \includegraphics[scale=0.3]{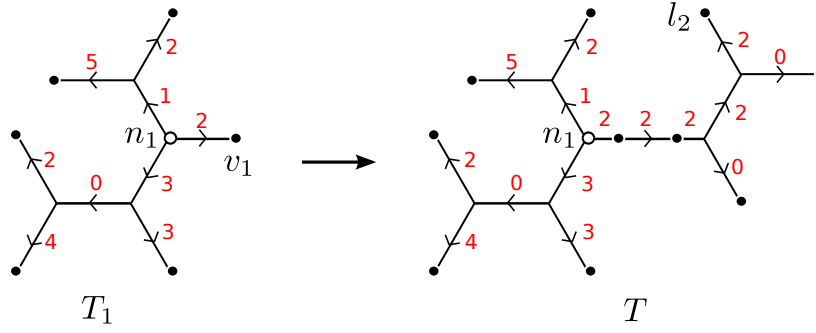}
 \label{extension}
\caption{A group-based flow for $\z_6$ on $T_1$ is extended to a group-based flow on $T$. }
\end{center}
\end{figure}

\begin{exm}
\rm
The figure \ref{extension} illustrates the above definition in the case $G=\z_6$ with an example of a flow in $T_1$ extended to a  flow in $T=T_1\star T_2$.
\end{exm}

Next, we proceed to define three sets of phylogenetic invariants for $T$.
Let $\AA_1$ (resp. $\AA_2$) be the set of phylogenetic invariants defining the variety $X_{T_1}$ (resp. $X_{T_2}$) on the respective Zariski open set $U_1:=U_{T_1}$ (resp. $U_2:=U_{T_2}$) as a complete intersection. In particular, $|\AA_i|=\g^{\ll_i-1}-1-(\g-1)\e_i$ (see Corollary \ref{cor:codimension}).

\begin{figure}
\begin{center}
 \includegraphics[scale=0.35]{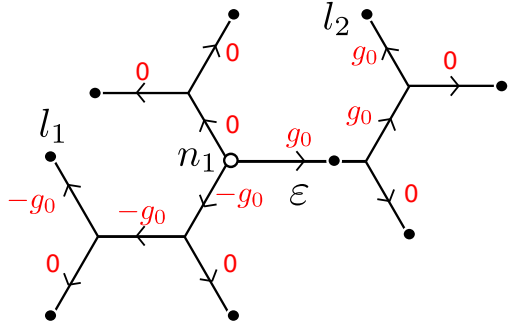}
\caption{\label{sp_flow} The group-based flow $\mathtt{f}_{g_0}$ defined on $T=T_1 \star T_2$.}
\end{center}
\end{figure}

\begin{enumerate}
\item[(1)] Invariants induced from $\AA_1$:  each invariant in $\AA_1$ is represented by two multisets of flows on $T_1$. Let us apply $E_1(\cdot)$ to all elements of both multisets, obtaining multisets $m_1$ and $m_2$ of flows on $T$. We claim that $m_1\equiv_T m_2$. This is equivalent to the equalities $\pi_e(m_1)=\pi_e(m_2)$ for every edge of $T$. This is obvious for edges in $T_1$, as before the extension we started from a valid relation on $T_1$. Notice that the value of the extension $E_1(f)$ on any edge of $T_2$ is uniquely determined by the element that $f$ associates to the pendant edge of $T_1$ where $v_1$ lies. As the projection to the chosen leaf gives the same multisets, the same must be true for all other edges of $T_2$. Hence, for each element of $\AA_1$ the multisets $m_1$ and $m_2$ define a phylogenetic invariant for $T$. Its degree is the same as the degree of the original element of $\AA_1$.

\vspace{2mm}
\item[(2)] Invariants induced from $\AA_2$: the construction is analogous to the previous case, by applying $E_2(\cdot)$.

\vspace{2mm}
\item[(3)] $\g(\g^{\ll_1-2}-1)(\g^{\ll_2-2}-1)$ quadratic invariants, which are examples of the so-called ``edge invariants" (cf. Example \ref{edgeinvariants}). These will come in $\g$ groups indexed by elements of $G$. Given $g_0\in G$, consider any group-based flow $f$ on $T$ that associates:
    \begin{itemize}
    \item $g_0$ to the common edge of $T_1$ and $T_2$ (there are $\g$ choices for these);
    \item  arbitrary elements to leaves of $T_1$ different from $v_1$, but not the neutral element at the same time to all leaves different from $l_1$. There are $\g^{\ll_1-2}-1$ possible choices for these.
    \item  arbitrary elements to leaves of $T_2$ different from $v_2$, but not the neutral element at the same time to all leaves different from $l_2$. There are $\g^{\ll_2-2}-1$ possible choices for these.
    \end{itemize}
There are $\g (\g^{\ll_1-2}-1)(\g^{\ll_2-2}-1)$ choices for $f$. Write $f_{|T_i}$ for the restriction of $f$ to $T_i$. Hence, $f$ can be considered as the join of $f_{|T_1}$ and $f_{|T_2}$: $f=\join{f_{|T_1}}{f_{|T_2}}$. %

\vspace{2mm}
\noindent
\textbf{Notation. }
We will write $\mathtt{f}_{g_0}$ for the group-based flow on $T=T_1\star T_2$  that assigns $g_0$ to all edges in the shortest path joining $n_1$ and $l_2$, $-g_0$ to all edges in the shortest path joining $n_1$ and $l_1$, and the neutral element to the other edges (see Figure \ref{sp_flow}).
Notice that for any flow $f$ as above, there is a quadratic relation:
\begin{eqnarray*}
\{f,\mathtt{f}_{g_0}\}\equiv \left \{ \join{f_{|T_1}}{{\mathtt{f}_{g_0}}_{|T_2}}, \join{{\mathtt{f}_{g_0}}_{|T_1}}{f_{|T_2}} \right \}
\end{eqnarray*}
where on the right hand side we have joins of respective restrictions.
\end{enumerate}

In summary,
we have defined
\begin{eqnarray*}
\left ( \g^{\ll_1-1}-1-(\g-1)e_1\right )+\left ( \g^{\ll_2-1}-1-(\g-1)e_2\right )+\g(\g^{\ll_1-2}-1)(\g^{\ll_2-2}-1)=\codim X
\end{eqnarray*}
invariants.
\lem
The invariants above form a set of Laurent monomials that define $X$ on $U$.
\klem
\dow
Consider any Laurent monomial 
vanishing on $X\cap U$, represented by multisets $m_1=\{f_1,\dots,f_k\}$ and $m_2=\{f'_1,\dots,f'_k\}$ of group-based flows on $T$ (that is, $\{f_1(e),\dots,f_k(e)\}=\{f'_1(e),\dots,f'_k(e)\}$  as multisets for any edge $e$ of $T$).
Consider the multisets
\begin{eqnarray*}
m'_1:=\{f_1,\dots,f_k,\mathtt{f}_{f_1(\varepsilon)},\dots,\mathtt{f}_{f_k(\varepsilon)} \},\\
m'_2:=\{f_1',\dots,f_k',\mathtt{f}_{f_1'(\varepsilon)},\dots,\mathtt{f}_{f_k'(\varepsilon)}\},
\end{eqnarray*}
where $\mathtt{f}_g$ is defined as above. %
Notice, that as $\pi_{\varepsilon}(m_1)=\pi_{\varepsilon}(m_2)$ (because $m_1$ and $m_2$ represent a Laurent monomial defining $X$ in $U$), we have enlarged $m_1$ and $m_2$ by adding the same multiset $\{\mathtt{f}_{f_1(\varepsilon)},\dots,\mathtt{f}_{f_k(\varepsilon)}\}$ of flows. Thus, as we consider the variety $X$ on the Zariski open set $U$, it is enough to see that the relation $m_1'\equiv m_2'$ can be generated by a relation in the elements of (1), (2), and (3). We can apply quadric relations
\begin{eqnarray*}
\{f_j,\mathtt{f}_{f_j(\varepsilon)}\}\equiv \left \{ \join{f_{j|T_1}}{\mathtt{f}_{f_j(\varepsilon)|T_2}},\join{\mathtt{f}_{f_j(\varepsilon)|T_1}}{f_{j|T_2}}\right \}
\end{eqnarray*}
and
\begin{eqnarray*}
\{f_j',\mathtt{f}_{f_j'(\varepsilon)}\}=\left \{\join{f_{j|T_1}'}{\mathtt{f}_{f_j'(\varepsilon)|T_2}},\join{\mathtt{f}_{f_j'(\varepsilon)|T_1}}{f_{j|T_2}'}\right \}
\end{eqnarray*}
for $j=1,\dots, k$.
After this reduction our relation $m_1'\equiv  m_2'$ is a sum of two relations:
\begin{eqnarray*}
\left \{ \join{f_{j|T_1}}{\mathtt{f}_{f_j(\varepsilon)|T_2}}\right \}_{j=1,\ldots,k}=\left \{\join{f_{j|T_1}'}{\mathtt{f}_{f_j'(\varepsilon)|T_2}}\right \}_{j=1,\ldots,k}
\end{eqnarray*}
and
\begin{eqnarray*}
\left \{\join{\mathtt{f}_{f_j(\varepsilon)|T_1}}{f_{j|T_2}}\right \}_{j=1,\ldots,k}=\left \{\join{\mathtt{f}_{f_j(\varepsilon)|T_1}}{f_{j|T_2}}\right \}_{j=1,\ldots,k}.
 \end{eqnarray*}

The first (resp. second) one is the extension $E_1(\cdot)$ (resp. $E_2(\cdot)$) of a relation holding on $T_1$ (resp. $T_2$). Hence, it is generated by binomials defined in point (1) (resp. (2)).
\kdow

\section{Complete intersection for claw trees}\label{clawtrees}
The varieties associated to trees of high valency are considered to be much more complex than those associated to trivalent trees. In this section we prove the following result, which gives a positive answer to a conjecture in the third author's PhD thesis \cite{jaPhD}.
\tw
The variety $X$ associated to the claw tree $T$ with $\ll$ leaves is a complete intersection in the Zariski open set $U$. Moreover, if $\ll\geq 4$, then $X$ is the scheme theoretic intersection in $U$ of two varieties associated to two trees of smaller valency, and we provide an explicit description of $X$ as a complete intersection in $U$.
\ktw

In order to prove this theorem, we shall consider two types of invariants. First, let $l_1$ and $l_2$ be two leaves of $T$ and consider a tree $T'$ with the same leaves as $T$ but with two interior  nodes: one of valency $3$ leading to leaves $l_1$ and $l_2$ and the other of valency $\ll-1$ leading to the rest of the leaves . Then the variety $X'$ associated to $T'$ contains $X$, so its defining equations on $U$ are also equations for $X$ (cf. Figure \ref{fig:2trees} and Example \ref{subtree}).

We now define $\g-1$ additional  phylogenetic invariants for the claw tree $T$, which will be quadrics indexed by non-neutral group elements. We consider $T$ rooted at the interior node.  We choose two leaves $l_3$, $l_4$ in $T$ different from $l_1$, $l_2$ and, without loss of generality, we assume that $l_1,\dots, l_4$ are the first four leaves in $T$.

As $T$ only contains pendant edges $e_1,\dots, e_{\ll}$, a group-based flow $f$ on $T$ will be denoted as $[g_1,\dots,g_{\ll}]$ if $f(e_i)=g_i$.  Note that the tuple $[g_1,\dots,g_{\ll}]$ of elements of $G$ is a group-based flow on $T$ if and only if $g_1+\dots+g_{\ll}=0$.  Suppose $G=\z_{a_1}\times\dots\times\z_{a_k}$, and
write $\m_j$ for the image of $m\in \z_{a_j}$ by the embedding $\z_{a_j} \hookrightarrow G$, that is, the element in $G$ whose entry in the $j$-th position is $m$ and the rest of entries are $0$.
%
%

Next, we proceed to assign a phylogenetic invariant to every element $b=
(b_1,\dots,b_k)\neq 0$ of $G$.
%


 \vspace{2mm}

\noindent
\emph{(Nonspecial).}
Assume $b$ is different of $\1_i$ for any $i=1,\dots,k$.
 Let $j$ be the largest index such that $b_j\neq 0$.  Define a quadratic relation $q_b$ of group-based flows as follows:
 \begin{eqnarray*}
  q_b: \{ lq^1_b,lq^2_b\} \equiv \{rq^1_b, rq^2_b \},
 \end{eqnarray*}
where
 \begin{eqnarray*}
 \begin{tabular}{lll}
  $lq_b^1 = [\1_j,0,b-\1_j,-b,0,\dots,0] $, & \qquad & $rq_b^1 = [0,0,b-\1_j,-b+\1_j,0,\dots,0], $\\
 $lq_b^2 = [0,b-\1_j,0,-b+\1_j,0,\dots,0]$, & \qquad & $rq_b^2 = [\1_j,b-\1_j,0,-b,0,\dots,0]$.
 \end{tabular}
  \end{eqnarray*}

%
%
%
\noindent
\emph{(Special).}
If $b=\1_j$ for some $j$, consider 
 \begin{eqnarray*}
  q_j: \{ lq^1_j,lq^2_j\} \equiv \{rq^1_j, rq^2_j \},
 \end{eqnarray*}
 where
  \begin{eqnarray*}
 \begin{tabular}{lll}
  $lq_j^1 = [\1_j,0,-\1_j,0,0,\dots,0] $, & \qquad & $rq_j^1 = [0,0,-\1_j,\1_j,0,\dots,0], $\\
 $lq_j^2 = [0,-\1_j,0,\1_j,0,\dots,0]$, & \qquad & $rq_j^2 = [\1_j,-\1_j,0,0,0,\dots,0]$.
 \end{tabular}
  \end{eqnarray*}
  In any case, these correspond to phylogenetic invariants on $T$ because the \emph{left} flows $lq$'s assign at each edge of $T$ that same pair as the \emph{right} flows $rq$'s. %
The last quadrics $q_j$, $j=1,\ldots,k$ will be called \emph{special} whereas the previous will be called \textit{nonspecial}.

Note that all these quadrics are edge-invariants for the tree $T''$ that has two interior nodes: one $w$ with descendants $l_2$ and $l_3$ and another with the rest of the leaves as descendants (cf. Example \ref{subtree}). Indeed, rooting $T''$ at $w$ for example,  all flows above can be extended to the internal edge of $T''$ by the same element: $\1_j$ in the special and $-b+\1_j$ in the nonspecial case.

\begin{figure}
\begin{center}
 \includegraphics[scale=0.35]{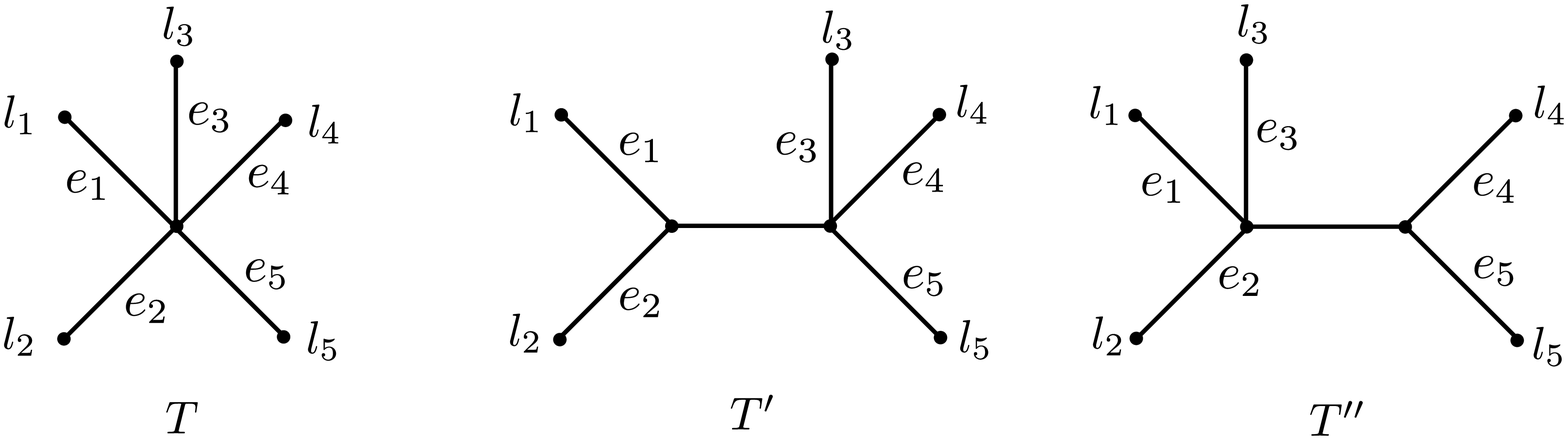}
\caption{\label{fig: high_valency} The trees $T$, $T'$ and $T''$ in the five leaves case.}
\end{center}
\end{figure}

\begin{exm}
If $T$ is the 5-leaved tree of fig \ref{fig: high_valency} and $G=\z_2\times \z_3$, then the quadric $q_b$ for the element $b=(1,2)$ is
\begin{eqnarray*}
 \big \{ [(0,1),(0,0),(1,1),(-1,-2),(0,0)],[(0,0),(1,1),(0,0),(-1,-1),(0,0)] \big \} \equiv \\
 \equiv \big \{[(0,0),(0,0),(1,1),(-1,-1),(0,0)],[(0,1),(1,1),(0,0),(-1,-2),(0,0)]\big \}.
\end{eqnarray*}
As all flows have $\pi_{e_2}+\pi_{e_3}=(1,1)$, this quadric is also a quadratic relation in ~$T''$.

\end{exm}
\dow
We proceed by induction on the number of leaves $\ll$. The case $\ll=3$ has been studied as a separate case in section \ref{sec:tripod}, so we assume $\ll>3$.

By the induction hypothesis, the variety associated to $T'$ on $U$ is a complete intersection defined by $\g^{\ll-1}-1-(\ll+1)(\g-1)$ phylogenetic invariants. All these are also invariants for the claw tree $T$, because $X_T$ is contained in $X_{T'}$.

These invariants together with the quadrics $q_b$, $b\neq 0$, defined above form a set of $\g^{\ll-1}-1-(\ll+1)(\g-1)+(\g-1)=\g^{\ll-1}-1-\ll(\g-1)=\codim X$ invariants. It remains to prove that on $U$ they generate any binomial in the ideal of $X$.  

Let us represent any such binomial by $m_1\equiv_T m_2$, where $m_1$ and $m_2$ are multisets of flows on $T$. The fact that it vanishes on $X$ is equivalent to the condition that the projection $\pi_{e_k}$ to any  edge $e_k$ of $T$  applied to $m_1$ and $m_2$ gives the same multisets of elements of $G$. Consider an operator $\pi_{1,2}$ that associates to any flow $f=[g_1,\dots,g_{\ll}]$ the sum $\pi_{1,2}(f)=g_1+g_2\in G$, 
and then represents it as an element $(b_1,\dots,b_k)\in \Z^k$, with $0\leq b_i<a_i$. We extend this operator $\pi_{1,2}$ to multisets of group-based flows as follows: for a multiset $m=\{f_1, \dots, f_d\}$ we define $\pi_{1,2}(m)$ to be the multiset $\{\pi_{1,2}(f_1),\dots,\pi_{1,2}(f_d)\}$ of elements in $\Z^k$.

Notice that if $m_1,m_2$ are multisets of flows on $T$ such that $\pi_{1,2}(m_1)=\pi_{1,2}(m_2)$,
then the binomial represented by $m_1$ and $m_2$ vanishes on $X_{T'}$: $m_1\equiv_{T'} m_2$, and hence it can be generated by the elements of a complete intersection for $T'$. %
In case $\pi_{1,2}(m_1)\neq \pi_{1,2}(m_2)$, we will reduce the multisets $m_1$ and $m_2$ using the $(\g-1)$ quadrics $q_b$ defined above, until $\pi_{1,2}$ applied to both multisets gives the same result.

If $\alpha$ is a multiset of elements in $\Z^k$, we denote by $s(\alpha)$ the sum of its elements, $s(\alpha)\in \Z^k$.

\vspace{2mm}

\noindent
\textsc{1st step.} We first show that any binomial represented by two multisets $m_1$ and $m_2$ can be replaced by a new binomial represented by multisets $m_1'$ and $m_2'$ that satisfy $s(\pi_{1,2}(m_1'))=s(\pi_{1,2}(m_2'))$ in $\z^k$.

We observe that although $\pi_{1,2}(m_1)$ and $\pi_{1,2}(m_2)$ may be different multisets, for sure one has $s(\pi_{1,2}(m_1))=s(\pi_{1,2}(m_2))$ as elements of $G$ (as $\pi_{e_i}(m_1)=\pi_{e_i}(m_2)$ as elements of $G$ for $i=1,2$). 
Therefore, although the sum of $\pi_{1,2}(m_1)$ and of $\pi_{1,2}(m_2)$ may not be the same vectors in $\z^k$, their difference in the $j$-th coordinates will be always divisible by $a_j$.

Note that the multisets $lq_j$ and $rq_j$ defining the special quadrics $q_j$ above satisfy $s(\pi_{1,2}(lq_j))=s(\pi_{1,2}(rq_j))+\aj_j$ (since $\pi_{1,2}(lq_{j})=\{\1_j,\aj_j-\1_j\}$ and $\pi_{1,2}(lq_{j})=\{0,0\}$).
Hence, by enlarging multisets $m_1$ and $m_2$ with the multisets $\{lq^1_j, lq^2_j\}$ and $\{rq^1_j,rq^2_j\}$ defined above respectively, we can assume that $\pi_{1,2}(m_1)$ and $\pi_{1,2}(m_2)$ sum up to the same vector in $\z^k$.

\vspace{2mm}

\noindent
\textsc{2nd step.} Now we assume that $s(\pi_{1,2}(m_1))=s(\pi_{1,2}(m_2))$ in $\z^k$ and we prove that $m_1,m_2$ can be replaced by two new multisets satisfying $\pi_{1,2}(m_1')=\pi_{1,2}(m_2')$.

To this end, we will use the nonspecial quadrics defined above. If $f$ is an element in $m_1$ such that $\pi_{1,2}(f)=(b_1,\dots,b_k)$ is different from zero or any $\1_i$, $i=1,\dots,k$, we define new multisets $m'_1= m_1\cup lq_b$ and $m'_2 = m_2\cup rq_b$, where
\begin{eqnarray*}
lq_b & = & \{lq^1_b, lq^2_b\}, \\ rq_b & = & \{rq^1_b, rq^2_b\}.
\end{eqnarray*}
 We have that $\pi_{1,2}(lq_b)=\{\1_j,b-\1_j\}$ and $\pi_{1,2}(rq_b)=\{0,b\}$.
In this case we say that $f$~and $rq_b^2$ are $\pi_{1,2}$-\emph{paired}. The other flows that have been added, $lq_b^1$, $lq_b^2$ and $rq_b^1$, are not $\pi_{1,2}$-paired, but their $\pi_{1,2}$ value is either $b-\1_j$, $0$, or $\1_j$. In any case, the corresponding $\pi_{1,2}$ value is smaller than $\pi_{1,2}(f)$. Moreover, $m_1'$ and $m_2'$ still satisfy $s(\pi_{1,2}(m_1'))=s(\pi_{1,2}(m_2'))$ in $\z^k$ because the multisets that have been added to $m_1$ and $m_2$ fulfill this condition  also.

We repeat the procedure, dealing also with the flows in $m_2$. In the end, we reach a couple of multisets $m'_1$ and $m'_2$, where the only (possibly) elements that are not $\pi_{1,2}$-paired elements are either equal to $0$ or $\1_j$ for some $i=1,\dots,k$. As  the sums of elements of $\pi_{1,2}(m_1')$ and $\pi_{1,2}(m_2')$ are equal as elements of $\z^k$, 
we deduce that $\pi_{1,2}(m_1')$ and $\pi_{1,2}(m_2')$ contain the same number of elements of type $\1_i$, the same number of elements equal to $0$, and a certain number of $\pi_{1,2}$-paired elements.  This means that we obtained a pair of multisets $m_1'$, $m_2'$ for which $\pi_{1,2}(m_1')=\pi_{1,2}(m_2')$ as multisets, as desired.

Such a relation is induced by a relation holding on $T'$, so we are done.
\kdow

\section{Conclusion}

Putting together all the results we have obtained,

\tw\label{th:main}
For any abelian group $G=\z_{a_1}\times\dots\times\z_{a_k}$ and any tree $T$ the associated variety $X$ in the Zariski open set $U$ is a complete intersection of explicitly constructed phylogenetic invariants of degree at most $\max(3,a_i)$.
\ktw

Varieties $X$ representing group-based models are complicated from an algebraic point of view. For arbitrary finite abelian group $G$ a complete description of the ideal is not known, even for the simplest tree (the tripod). Also, for simple groups, like $\z_2\times\z_2$ a complete description of the ideal is only conjectural for arbitrary trees. However, these varieties admit a simple description on the Zariski open set $U$, isomorphic to a torus. In Fourier coordinates this torus is identified with the locus of points in the projective space with all coordinates different from zero. The intersection $X\cap U$ is a torus, which admits a precise description as a complete intersection (in $U$) of $\codim X$ phylogenetic invariants of degree at most $|G|$. Thus for a fixed $G$, to find these phylogenetic invariants explicitly one proceeds as follows:
\begin{enumerate}
\item present $G$ as a product of cyclic groups,
\item for each cyclic component of $G$ find the correct phylogenetic invariants for the tripod -- the explicit formula for them is given in the proof of Proposition \ref{prop:cyclic},
\item reconstruct the correct phylogenetic invariants for the tripod and the whole group $G$ -- these amounts to adding specific cubics, as described in Section \ref{arbitgr},
\item if we consider any trivalent tree an inductive procedure, basing on adding correct edge invariants, to construct correct phylogenetic invariants was given in Section \ref{gluetrees},
\item if we want to find phylogenetic invariants for trees of higher valency, we have to first construct them for claw trees (the method is provided in Section \ref{clawtrees}) and as before apply results of Section \ref{gluetrees}.
\end{enumerate}
In particular, on $U$, the conjecture of Sturmfels and Sullivant on the degree of phylogenetic invariants holds.

\section{Appendix}
\dow[Proof of Proposition \ref{cor:codimension}]
The last part of the Proposition is implied by:
$$\c[\tilde M_0]=\c[M_0]^{(G^N)}.$$
Thus it is enough to prove the above equality.

Clearly the elements of $\tilde M_0$ are invariant under the action of $G^N$, hence $\c[\tilde M_0]\subset \c[M_0]^{(G^N)}$. The elements of $M_0$ form a basis of $\c[M_0]$ consisting of eigenvectors with respect to the $G^N$ action. Thus any invariant vector must be a linear combination of invariant elements of $M_0$. It remains to prove that an element of $M_0$ that is invariant with respect to $G^N$ belongs to $\tilde M_0$. The proof is inductive on the number of nodes of the tree $T$.

First suppose that $T$ has one interior node, that is $T$ is a claw tree, with $\ll$ leaves. Consider an invariant element of $M_0$ given by $R:=\sum_{j=1}^\ll\sum_{g\in G} a_{(j,g)} b_{(j,g)}$ with the condition $\sum_{g\in G} a_{(1,g)}=\dots=\sum_{g\in G} a_{(\ll,g)}=0$. We will reduce $Q$ to zero modulo $\tilde M_0$. Notice that for any $1\leq j\leq\ll$, $g_1,g_2\in G$ the element $S_{j,g_1,g_2}:=b_{(j,g_1)}+b_{(j,g_2)}-b_{(j,g_1+g_2)}-b_{(j,\0)}$ belongs to $\tilde M_0$. Indeed, for example for $j=1$ it equals:
$$Q_{[g_1,-g_1,0,\dots,0]}+Q_{[g_2,0,-g_2,0,\dots,0]}-Q_{[g_1+g_2,-g_1,-g_2,0,\dots,0]}-Q_{[0,\dots,0]}.$$

Using elements as above we can reduce $R$ and assume that for any $g\neq 0$ and $1\leq j\leq \ll$, the coefficient $a_{(j,g)}$ is zero apart from one $g$ for each $j$, for which the coefficient can be equal to one. Precisely, if for some $j$ coefficients $a_{(j,g_1)},a_{(j,g_2)}$ are positive (resp. negative) we subtract (resp. add) $S_{j,g_1,g_2}$. If there is a positive entry $a_{(j,g_1)}$ and a negative $a_{(j,g_2)}$ we add $S_{j,g_2,g_1-g_2}$. If a coefficient $a_{(j,g)}$ is negative we add $S_{j,g,-g}$. If a coefficient $a_{j,g}>1$ we subtract $S_{j,g_1,g_1}$. All these operations either strictly decrease $\sum_{g\neq0} |a_{j,g}|$ or leave the sum unchanged and increase the sum of negative coefficients. Thus the procedure must finish.

In other words, $R=\sum_{j=1}^\ll b_{(j,g_j)}-Q_{[0,\dots,0]}$ modulo $\tilde M_0$. As $R$ is invariant, we obtain $\sum_{j=1}^\ll g_j=0$, which finishes the first inductive step.

Suppose now that $T$ has more than one interior nodes. Consider an invariant element $R\in M_0$ as before. By choosing an interior edge $m\in E$ we can present $T=T_1\star T_2$. The element $Q$ induces two invariant elements $R_i\in M_{0,T_i}$ for $i=1,2$. By the inductive assumption we obtain: $R_i=\sum_j c_{i,j}Q_{f_{i,j}}$, where $c_{i,j}\in\z$, $\sum_j c_{i,j}=0$ and $Q_{f_{i,j}}\in P_{T_i}$ correspond to flows $f_{i,j}$ on the tree $T_i$. Let us consider the signed multisets\footnote{Formally, by a signed multiset we mean a pair of multisets on the same base set. The first multiset represents the positive multiplicities, the second one negative.} $Z_i$ that are the projections of $\sum c_{i,j}Q_{f_{i,j}}$ onto the edge $m$ -- each $f_{i,j}$ distinguishes an element on $m$. The multiset $Z_i$ has $c_{i,j}$ elements distinguished by $f_{i,j}$ with a minus sign if $c_{i,j}<0$. $Z_i$ is a signed multiset of group elements. Let $Z_i'$ be a signed multiset obtained by reductions cancelling $g$ with $-g$ in the multiset $Z_i$\footnote{Formally, if an element belongs to both multisets (the negative and the positive one) we cancel it.}. The multiset $Z_1'$ is just the signed multiset of group elements corresponding to the projection of $R$ to $m$. Thus, the multiset $Z_2'$ is the same multiset as $Z_1'$. This means that we can pair together elements from $Z_1'$ and $Z_2'$ such that each pair gives rise to a flow on the tree $T$. The image of the sum of these flows does \emph{not} have to equal $R$ yet. We have to lift also the flows that we canceled by passing from $Z_i$ to $Z_i'$. This is done as follows. Suppose that two flows $f_{1,j_0}$ and $f_{1,j_1}$ on $T_1$ associate $g$ to the edge $m$, but $c_{1,j_0}>0$ and $c_{1,j_1}<0$. Then, $f_{1,j_0}$ and $-f_{1,j_1}$ were canceling each other in $Z_1$. We choose any flow $s$ on $T_2$ that associates $g$ to the edge $m$. We can glue together $f_{1,j_0}$ and $s$ obtaining a flow $f_{1,j_0}\star s$ on the tree $T$ and analogously $f_{1,j_1}\star s$. The difference of flows $Q_{f_{1,j_0}\star s}-Q_{f_{1,j_1}\star s}$ has the same coordinates $b_{(e,g)}$ on the edges $e$ of the tree $T_1$ as $Q_{f_{1,j_0}}-Q_{f_{1,j_1}}$. Moreover, the coordinates $b_{(e,g)}$ for the edges $e$ belonging to $T_2$ are equal to zero. In this way we obtain the flows of $T$ with the signed sum equal to $\sum c_jf_{i,j}$ on $T_i$, hence equal to $R$.
\kdow
\bibliographystyle{amsalpha}
\bibliography{Xbib}

\def\cprime{$'$}
\providecommand{\bysame}{\leavevmode\hbox to3em{\hrulefill}\thinspace}
\providecommand{\MR}{\relax\ifhmode\unskip\space\fi MR }
\providecommand{\MRhref}[2]{%
  \href{http://www.ams.org/mathscinet-getitem?mr=#1}{#2}
}
\providecommand{\href}[2]{#2}
\begin{thebibliography}{CFS11}

\bibitem[AR03]{Allman2003}
ES~Allman and JA~Rhodes, \emph{Phylogenetic invariants for the general {M}arkov
  model of sequence mutation}, Mathematical Biosciences \textbf{186} (2003),
  no.~2, 113--144.

\bibitem[AR04]{Allman2004}
\bysame, \emph{Quartets and parameter recovery for the general {M}arkov model
  of sequence mutation}, Applied Mathematics Research Express \textbf{2004}
  (2004), no.~4, 107--131.

\bibitem[AR07]{AllmanRhodes_chapter4}
E~S Allman and J~A Rhodes, \emph{Phylogenetic invariants}, Reconstructing
  Evolution (O~Gascuel and MA~Steel, eds.), Oxford University Press, 2007.

\bibitem[AR08]{AllRhMarkov}
Elizabeth~S. Allman and John~A. Rhodes, \emph{Phylogenetic ideals and varieties
  for the general {Markov} model}, Advances in Applied Mathematics
  \textbf{40(2)} (2008), 127--148.

\bibitem[BW07]{BW}
Weronika Buczy\'nska and Jaros{\l}aw~A. Wi\'{s}niewski, \emph{On geometry of
  binary symmetric models of phylogenetic trees}, J. Eur. Math. Soc.
  \textbf{9(3)} (2007), 609--635.

\bibitem[Cas12]{CasanellasNewsEMS}
M~Casanellas, \emph{Algebraic tools for evolutionary biology}, EMS Newsletter
  \textbf{86} (2012), 12--18.

\bibitem[CFS07]{CFS_MBE}
M~Casanellas and J~Fernandez-Sanchez, \emph{Performance of a new invariants
  method on homogeneous and nonhomogeneous quartet trees}, Mol. Biol. Evol.
  \textbf{24} (2007), no.~1, 288--293.

\bibitem[CFS08]{CFS}
M.~Casanellas and J.~Fernandez-Sanchez, \emph{Geometry of the {Kimura}
  3-parameter model}, Advances in Applied Mathematics \textbf{41} (2008),
  265--292.

\bibitem[CFS11]{CFS3}
M~Casanellas and J~Fernandez-Sanchez, \emph{Relevant phylogenetic invariants of
  evolutionary models}, Journal de Mathématiques Pures et Appliquées
  \textbf{96} (2011), 207--229.

\bibitem[CLS11]{Cox}
David~A Cox, John~B Little, and Henry~K Schenck, \emph{Toric varieties},
  American Mathematical Soc., 2011.

\bibitem[Coh04]{Cohen}
Joel~E Cohen, \emph{Mathematics is biology's next microscope, only better;
  biology is mathematics' next physics, only better}, PLoS Biol \textbf{2}
  (2004), no.~12.

\bibitem[CP07]{Sonja}
J.~Chifman and S.~Petrovi\'{c}, \emph{Toric ideals of phylogenetic invariants
  for the general group-based model on claw trees $k_{1,n}$}, Proceedings of
  the 2nd international conference on {Algebraic} biology (2007), 307--321.

\bibitem[DBM]{DBM}
Maria Donten-Bury and Mateusz Micha{\l}ek, \emph{Phylogenetic invariants for
  group-based models}, arXiv:1011.3236v1.

\bibitem[DK09]{DK}
Jan Draisma and Jochen Kuttler, \emph{On the ideals of equivariant tree
  models}, Mathematische Annalen \textbf{344(3)} (2009), 619--644.

\bibitem[Ful93]{Fult}
William Fulton, \emph{Introduction to toric varieties}, Annals of Mathematics
  Studies, vol. 131, Princeton University Press, Princeton, NJ, 1993, The
  William H. Roever Lectures in Geometry.

\bibitem[HP89]{HendyPenny}
Michael Hendy and David Penny, \emph{A framework for the quantitative study of
  evolutionary trees}, Systematic Zoology \textbf{38} (1989), 297--309.

\bibitem[Mic11]{JaJalg}
Mateusz Micha{\l}ek, \emph{Geometry of phylogenetic group-based models},
  Journal of Algebra \textbf{339} (2011), no.~1, 339--356.

\bibitem[Mic12]{jaPhD}
\bysame, \emph{Toric varieties: phylogenetics and derived categories}, PhD
  thesis (2012).

\bibitem[Mic13]{JaJCTA}
\bysame, \emph{Constructive degree bounds for group-based models}, Journal of
  Combinatorial Theory, Series A \textbf{120} (2013), no.~7, 1672--1694.

\bibitem[Mic14]{JaAdvgeom}
\bysame, \emph{Toric geometry of the 3-kimura model for any tree}, Advances in
  Geometry \textbf{14} (2014), no.~1, 11--30.

\bibitem[PS04]{Pachter2004}
L~Pachter and B~Sturmfels, \emph{Tropical geometry of statistical models},
  Proceedings of the National Academy of Sciences \textbf{101} (2004),
  16132--16137.

\bibitem[PS05]{PS}
Lior Pachter and Bernd Sturmfels, \emph{Algebraic {Statistics} for
  {Computational} {Biology}}, Cambridge University Press, 2005.

\bibitem[SS05]{SS}
Bernd Sturmfels and Seth Sullivant, \emph{Toric ideals of phylogenetic
  invariants}, J. Comput. Biology \textbf{12} (2005), 204--228.

\bibitem[Stu96]{Stks}
Bernd Sturmfels, \emph{Gr\"obner bases and convex polytopes}, University
  Lecture Series, vol.~8, American Mathematical Society, 1996.

\bibitem[Sul07]{SethTFP}
Seth Sullivant, \emph{Toric fiber products}, Journal of Algebra \textbf{316}
  (2007), no.~2, 560 -- 577, Computational Algebra.

\end{thebibliography}
\end{document}